\documentclass[10pt]{article}
\usepackage{authblk}
\usepackage{graphicx}
\usepackage{epstopdf}

\begin{document}

\title{An observation on \textquotedblleft Classification of Lie point
symmetries for quadratic Li\'{e}nard type equation $\ddot{x}+f\left(
x\right) \dot{x}^{2}+g\left( x\right) =0$\textquotedblright\ [J. Math.
Phys. 54, 053506 (2013)] and its erratum [J. Math. Phys. 55, 059901 (2014)]}
\author[1]{A Paliathanasis\thanks{%
anpaliat@phys.uoa.gr}}
\author[2,3,4]{ PGL Leach\thanks{%
leach@ucy.ac.cy}}

\affil[1]{Instituto de Ciencias F\'{\i}sicas y Matem\'{a}ticas, Universidad Austral de Chile, Valdivia, Chile}
\affil[2]{Department of Mathematics and
Institute of Systems Science, Durban
University of Technology, PO Box 1334, Durban 4000, Republic of South Africa}
\affil[3]{School of Mathematics, Statistics and Computer Science, University
of KwaZulu-Natal, Private Bag X54001, Durban 4000, Republic of South Africa}
\affil[4]{Department of Mathematics and Statistics, University of Cyprus,
Lefkosia 1678, Cyprus}

\renewcommand\Authands{ and }

\maketitle
\begin{abstract}
We demonstrate a simplification of some recent works on the classification of
the Lie symmetries for a quadratic equation of Li\'{e}nard type. We observe
that the problem could have been resolved more simply.
\end{abstract}

\bigskip In \cite{ref1} and \cite{ref2} the classification of the Lie
(point) symmetries for the quadratic equation of the form%
\begin{equation}
\ddot{x}+f\left( x\right) \dot{x}^{2}+g\left( x\right) =0,  \label{eq1}
\end{equation}%
was performed, where overdot denotes total differentiation with respect to
time, \textquotedblleft $t$\textquotedblright .

We observe that under the coordinate (not point) transformation,%
\begin{equation}
y=\int \left[ \exp \left( \int f\left( x\right) dx\right) \right] dx,
\label{eq2}
\end{equation}%
equation (\ref{eq1}) takes the simpler form%
\begin{equation}
\ddot{y}+F\left( y\right) =0,  \label{eq3}
\end{equation}%
where%
\begin{equation}
g\left( x\right) =e^{-\int f\left( x\right) dx}F\left( \int e^{\int f\left(
x\right) dx}dx\right) .  \label{eq.4}
\end{equation}

The classification of the Lie symmetries of the latter equation has been
known for some time (approximately 140 years). The various possibilities are

\begin{enumerate}
\item $F\left( y\right) $ is an arbitrary function. Equation (\ref{eq3})
admits the autonomous symmetry, $\partial_t$, and the equation can be simply
reduced to a first integral which in general cannot be evaluated to obtain
the solutions in closed form.

\item In the two cases

\begin{enumerate}
\item $F\left( y\right) =\left( \alpha +\beta y\right) ^{n},~n\neq 0,1,-3$
and

\item $F\left( y\right) =e^{\gamma y},~\gamma \neq 0,$
\end{enumerate}

the admitted Lie symmetries of (\ref{eq3}) constitute the algebra $A_{2}$ in
the Mubarakzyanov Classification Scheme \cite%
{Morozov58a,Mubarakzyanov63a,Mubarakzyanov63b,Mubarakzyanov63c}.

\item When

\begin{enumerate}
\item $F\left(y\right) =\frac{1}{\left( y+c\right) ^{3}}$ or

\item $F\left( y\right) =\alpha\left( y+c\right) +\frac{\beta }{\left(
y+c\right) ^{3}}$, $\beta\neq 0$,
\end{enumerate}

equation (\ref{eq3}) is invariant under the three-dimensional algebra, $%
A_{3,8}$, which is more commonly known as $sl\left( 2,R\right) $.

\item Finally for the cases

\begin{enumerate}
\item $F\left( y\right) = 0$,

\item $F\left( y\right) = c$,

\item $F\left( y\right) = y$ and

\item $F\left( y\right) = y + c$
\end{enumerate}

(note that in (c) and (d) a multiplicative constant -- or arbitrary function
of time which is beyond the considerations of \cite{ref1, ref2} -- in the $y$
term is superfluous), the algebra of the Lie symmetries is $sl\left(
3,R\right) $ and, as a second-order linear equation, (\ref{eq3}) is
maximally symmetric \cite{Lie67a}[p 405].
\end{enumerate}

\subsubsection*{Acknowledgements}

AP acknowledges Prof. PGL\ Leach, Sivie Govinder, as also DUT for the
hospitality provided and the UKZN, South Africa, for financial support. The
research of AP was supported by FONDECYT postdoctoral grant no. 3160121.

\end{document}